# Optimal Adaptive Droop Design via a Modified Relaxation of the OPF

H. Sekhavatmanesh, *Member, IEEE,* G. Ferrari-Trecate, *Senior Member, IEEE,* and *S. Mastellone, Member, IEEE.*

*Abstract*— **The ever increasing penetration of Renewable Energy Resources (RESs) in power distribution networks has brought, among others, the challenge of maintaining the grid voltages within the secure region. Employing droop voltage regulators on the RES's inverters is an efficient and low cost solution to reach this objective. However, fixing droop parameters or optimizing them only for overvoltage conditions does not provide the required robustness and optimality under changing operating conditions. In this paper, a convex optimization approach is proposed for reconfiguring PV and QV droop regulators during online operation. The objective is to minimize power curtailment, power losses, and voltage deviation subject to electrical security constraints. This enables to optimally operate the grid with high RES penetration under variable conditions while preserving electrical security constraints (e.g., current and voltage limits). As a first contribution, a mixed integer linear model of the droop characteristics is developed. According to this model, a droop regulated generation unit is represented as a constant power generator in parallel with a constant impedance load. As a second contribution, a *Modified Augmented Relaxed Optimal Power Flow* (MAROPF) formulation is proposed to guarantee that the electrical security constraints are respected in the presence of constant impedance loads in the network. Sufficient conditions for the feasibility of the MAROPF solution are provided. Those conditions can be checked *a priori* and are valid for several real distribution networks. Furthermore, an iterative approach is proposed to derive an approximate solution to the MAROPF that is close to the global optimal one. The performance of the MAROPF approach and the accuracy of the proposed model are evaluated on standard 34 bus and 85 bus test networks.**

*Index Terms*— **constant-impedance loads, droop control, optimal power flow (OPF), radial power grids.[1]**

## I. Introduction

The quest to address the current global energy challenges has led to a paradigm shift in power generation with an increased penetration of RESs in power distribution networks. This transition requires to rethink several control aspects including voltage regulation. To this end, different operational solutions have been proposed by the power systems and control research communities, e.g., network reconfiguration and optimal setting of conventional voltage regulation devices such as tap changing transformers. Compared to these conventional solutions, the control of Inverter-Based Distributed Generations (IBDGs) offers the benefits of low investment costs and fast dynamic responses.

The voltage control of IBDGs is typically addressed via: I) centralized approaches, which assume the complete knowledge of the network [1], [2] or II) local approaches, which solely rely on local voltage and power measurements [3], [4]. In the latter category, IBDGs can promptly react to the rapid change in voltages, by adjusting real and reactive power injections. According to the IEEE 1547.8 standard [5], the power adjustments should be carried out according to the Volt/Var and Volt/Watt (Q-V and P-V) droop characteristics depicted in Figure. 1. The IEEE 1547.8 standard also recommends typical values for the parameters of the droop curves. However, using fixed parameters limits the effectiveness of the local droop controller in achieving global objectives related to the entire network operation such as minimization of power loss [6], and total operational cost [7].

To address these shortcomings, several methodologies have been developed in the literature to optimize the droop parameters. For example, the droop slope is optimized in [8] to achieve a desirable load sharing, while the authors of [9] set the droop reference to minimize the total power loss in the network. In [10], an optimal dispatch problem is solved to find the optimal reactive power setpoints. There, a piece-wise linear Q-V characteristic is fit to the obtained reactive power setpoints, using machine learning methods. However, there is no guarantee of the feasibility and optimality of the obtained solution. Another methodology is proposed in [11], where the model of the P-V droop characteristics is approximated, assuming a pure linear characteristic during the entire voltage range.

An additional aspect to consider is the network electrical security that requires the nodal voltage and line current magnitudes in the whole grid to be bounded within prescribed limits. To model the AC power flow in an electric grid, these variables are expressed as functions of the power injections at some of the buses. To deal with the nonlinearity and nonconvexity inherent in the AC power flow equations, there are currently four approaches: I) considering only the economic aspects and disregarding the electrical security of the network

The manuscript is submitted on September 2nd, 2021. This work was supported by the Swiss National Science Foundation (SNSF) under the NCCR Automation project (grant agreement: 51NF40_180545).

H. Sekhavatmanesh (corresponding author) and S. Mastellone are with Institute of Electrical Engineering, FachHochschule NordWestschweiz (FHNW), Windisch, Switzerland (emails: hossein.sekhavat@fhnw.ch; silvia.mastellone@fhnw.ch). G. Ferrari-Trecate is with Institute of Mechanical Engineering, École Polytechnique Fédérale de Lausanne (EPFL), Switzerland (email: giancarlo.ferraritrecate@epfl.ch).

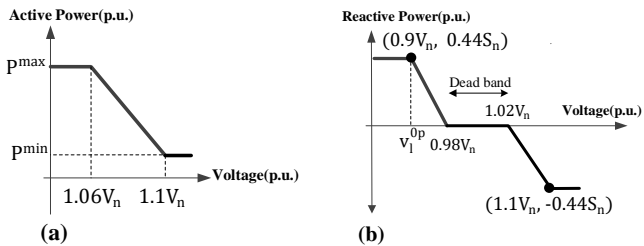

Figure. 1. The P-V (a) and Q-V (b) droop curves with the predefined parameters defined by IEEE 1547.8 standard

[12], II) using linearized formulation such as *DC power flow* [13], III) applying nonconvex formulation of the power flow, which are solved using either non-linear solvers [14] or heuristic-based solvers [15], and IV) applying convexification methods [16], which includes two classes. The first class is based on the convex relaxation and provides an outer approximation of the feasible set of the OPF problem. Several approaches have been proposed, e.g., chordal relaxation [17], semidefinite programming (SDP) [18], [19], least squares estimation-based SDP [20], and *Second-Order Cone Programming* (SOCP) [21], [22]. However, using only relaxation techniques in the optimization problem produces a lower bound of the objective function. This solution might violate the security constraints of the network, especially during high power production periods, where the upper voltage limit is binding.

To tackle this issue, several studies developed specific conditions that guarantee the exactness of the relaxed formulation. The authors of [21] introduced two conditions for the exactness of the SOCP relaxation, namely I) the voltage angle difference at the line's terminals is sufficiently small, and II) the nodal reactive power injection is not bounded. The latter assumption is usually not fulfilled in practice. In [22], it is proved that the SOCP relaxation is exact under a specific condition on the network parameters together with the condition that there is no upper bound on the voltage magnitude, which is hardly met in practice. As another example, the authors of [23] prove the exactness of the SOCP relaxation under a set of conditions. One of those requires that both active and reactive powers in each line do not flow in the reversed direction. Unfortunately, this condition is not always met, especially with the high penetration of distributed generation. According to [24], the optimal solution of the relaxed SOCP is exact if the objective function is strictly monotonic. This assumption does not always hold, e.g. in a case, where the objective is to track the reference of power injection from the upper grid [25]. To address this issue, the authors of [25] used a first-order approximation of the line losses formulation to make the cost function strictly monotonic. However, this approximation is not accurate and might lead to infeasible solutions. Moreover, the effect of reactive power dispatch and line current capacity was not considered.

The second class of convexification methods are restriction methods that aim to obtain an inner approximation of the feasibility set. Unlike relaxation techniques, these approaches provide sufficient conditions for guaranteeing the feasibility of the solution. Algorithms of this kind have been first developed in [26] for optimizing the reactive power setpoint of discrete reactive power controllers in radial and inductive networks. The objective was to maximize the voltage security margins. This approach has been extended in [27] so as to be applicable to any radial feeder with a mix of inductive and capacitive branches. The developed methodology is then used to derive nodal injection limits for each IBDG, such that the security constraints in the network are satisfied. In [28], [29], the authors developed a tractable convex restriction strategy for OPF problem in a general meshed network. All these restriction methodologies rely on the derivation of envelopes over the nonlinear terms of power flow equations. These envelopes are constructed around a fixed operating point of the network. The conservativeness of different methods depends on this operating point, and the obtained solution might be far from the global optimal one. To reduce the conservativeness of the solution, an iterative algorithm is proposed in [27], where this operating point is modified at each iteration. However, there is no guarantee that this method improves the quality of the solution for all objective functions and operating conditions.

For achieving both feasibility and optimality by using a convex optimization approach for the OPF problem, an idea is to employ both relaxation and restriction methodologies. This idea is employed in [24], augmenting the SOCP relaxed formulation with conservative bounds on the voltage and line current magnitudes. The exactness of the resulting Augmented Relaxed OPF (AR-OPF) formulation is proved under certain conditions that can be checked *a priori*. However, constant impedance loads in the network were not considered. The integration of ZIP load models into OPF formulations has been already studied in the literature, e.g., in [30], where constant impedance loads are integrated into SDP. However, to our knowledge, an analytical proof for the feasibility of the resulting optimization problem is not available in the literature.

To overcome some of the aforementioned limitations, in this paper, we propose a multi-period optimization-based approach to define the optimal droop parameters, where:

- the piecewise affine characteristic of the PV droop is modelled using mixed-integer linear constraints. According to this model, a droop-based generation unit is represented with a constant-power source in parallel with a constant-impedance load.
- to guarantee the electrical security constraints (e.g. nodal voltage and line current limits) a new formulation, called MAR-OPF, is proposed. This convex relaxation method can be applied in general cases, where upper-voltage and -current limits must be fulfilled in the presence of constant-impedance loads in the grid.
- the feasibility of the MAR-OPF solution is proved under some conditions that can be checked a priori.
- an iterative algorithm is proposed to improve the quality of the MAR-OPF solution.

The paper is organized as follows. In section II, the proposed MAR-OPF formulation for the convexification of the powerflow equations is presented. The developed MAR-OPF formulation is applied in section III to optimize the parameters of the P-V and Q-V droop controllers. To improve the quality of the obtained solution, an iterative approach is presented in section IV. Section V presents the results of the proposed optimization model applied to an IEEE standard test distribution network. Finally, concluding comments are provided in section VI.

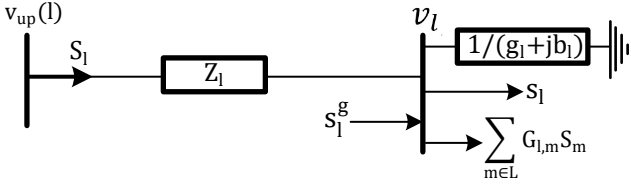

Figure. 2. The model of a distribution line with the notations used for the formulation of the power flow equations.

## II. Convex OPF formulation

### A. Notations

In the following, each distribution line is modeled with an impedance as shown in Figure. 2. The set of lines is represented with $\mathcal{L}$ and the index $l$ is used to refer both to a given line and to its ending node. The starting node of line $l$ is denoted with $up(l)$. The set of downstream buses of bus $l$ is $\mathcal{L}^l$ ($l \notin \mathcal{L}^l$).[1] The index 0 is reserved for the slack bus, whose voltage is fixed to $\sqrt{v_0}$. $\Re(.)$ represents the real part of a complex number.

$S_l = P_l + jQ_l$ denotes the complex power flow entering line $l$ from bus $up(l)$; $v_l$ and $f_l$ are the squares of the voltage and current magnitudes at bus $l$ and through line $l$, respectively; $v_l^{min}$ and $v_l^{max}$ are the square of the minimum and maximum voltage limits on bus $l$. $l_l^{max}$ denotes the square of the maximum current limit of line $l$. The maximum active and reactive power flows at each line $l$ are $P_l^{max}$ and $Q_l^{max}$, respectively. $p_{l,max}^g$ and $q_{l,max}^g$ denote the maximum active and reactive power limits of the IBDG at node $l$, respectively. $z_l = r_l + jx_l$ is the longitudinal impedance of line $l$. $g_l$ and $b_l$ represent the conductance and susceptance of the constant-impedance load at bus $l$. We denote by $s_l$ the complex power consumption at bus $l$, while the complex power injection is given by $s_l^g = p_l^g + jq_l^g$. The notations without subscript, such as $v$, denote the corresponding column vector, whose $l^{th}$ element is $v_l$.

In this paper bold nonitalic letters refer to the matrices, i.e., $\mathbf{I}$, which represents the identity matrix. The adjacency matrix of the directed graph of the distribution network is denoted with $\mathbf{G}$, meaning that $\mathbf{G}_{k,l}$ equals to one if $k = up(l)$, and equals to zero otherwise. Diagonal elements of $\mathbf{G}$ are zero. $\mathbf{H} = (\mathbf{I} - \mathbf{G})^{-1}$ is the closure of $\mathbf{G}$, i.e. $\mathbf{H}_{k,l}$ equals to one if bus $k$ is on the path from the slack bus to bus $l$ or $k = l$, and equals to zero otherwise. The Frobenius norm of a given matrix $\mathbf{A}$ is denoted with $\|\mathbf{A}\|$. diag(r) represents a diagonal matrix, whose $(l, l)$ element is $r_l$.

### B. General Relaxed OPF (R-OPF)

R-OPF refers to the optimization problem including the general SOCP relaxation of power flow equations, which are given below for each bus (or line) $l \in \mathcal{L}$ in a radial power grid. In the R-OPF problem, the unknown variables are $(s^g, S, v, f)$. Other symbols denote constant parameters.

$$S_l = s_l - s_l^g + v_l(g_l + jb_l) + \sum_{m \in \mathcal{L}} G_{l,m} S_m + Z_l f_l \quad (1.a)$$

$$v_l = v_{up(l)} - 2\Re(z_l^* S_l) + |z_l|^2 f_l \quad (1.b)$$

---

[1] The distribution network is assumed radial and balanced.

$$f_l = \frac{|S_l|^2}{v_l} \quad (1.c)$$

Equations (1.a) and (1.b) describe the Kirchhoff's current and voltage laws for each bus and line $l$, respectively. The current flow equation (1.c) is non-linear and makes the optimization problem non-convex. Therefore, in the R-OPF formulation [15], it is replaced with:

$$f_l \geq \frac{|S_l|^2}{v_l} \quad (2)$$

This relaxation renders R-OPF convex. However, where the maximum current or voltage limits are binding, the optimal solution of R-OPF might not satisfy the original constraint (1.c). In the following, we propose a MAR-OPF formulation, where additional constraints are added to avoid such inexact solutions, which have no physical correspondence in reality.

### C. Modified Augmented Relaxed OPF (MAR-OPF)

Inspired by the AR-OPF formulation proposed in [31], the main idea of MAR-OPF is to impose maximum voltage and current limits on a set of auxiliary variables ($\bar{f}, \hat{S}, \underline{S}, \bar{S}, \hat{v}$) to ensure the feasibility of the solution. Variables $\bar{f}$ and $\hat{v}$ represent upper bounds on-line current and nodal voltage magnitudes, which do not depend on the line current flow $f$. The formulation provided in [31] is extended here by accounting for the constant-impedance loads in the grid. $\hat{S} = \hat{P} + j\hat{Q}$ and $\hat{v}$ are defined in (3.a) and (3.b) according to the DistFlow equations, where the line current flows are assumed zero. The definition that is provided in [31] for auxiliary variables $\bar{S} = \bar{P} + j\bar{Q}$ and $\bar{f}$ is updated in (3.d) and (3.e) so that they are still upper bounds respectively for $S = P + jQ$ and $f$, even in presence of constant-impedance loads $g_l + jb_l$. In addition, new auxiliary variables $\underline{S} = \underline{P} + j\underline{Q}$ are added and defined in (3.c) as lower bounds for $S = P + jQ$.

$$\hat{S}_l = s_l - s_l^g + \hat{v}_l(g_l + jb_l) + \sum_{m \in \mathcal{L}} G_{l,m} \hat{S}_m \quad (3.a)$$

$$\hat{v}_l = \hat{v}_{up(l)} - 2\Re(z_l^* \hat{S}_l) \quad (3.b)$$

$$\underline{S}_l = s_l - s_l^g + v_l(g_l + jb_l) + \sum_{m \in \mathcal{L}} G_{l,m} \underline{S}_m \quad (3.c)$$

$$\bar{S}_l = s_l - s_l^g + v_l(g_l + jb_l) + \sum_{m \in \mathcal{L}} G_{l,m} \bar{S}_m + Z_l \bar{f}_l \quad (3.d)$$

$$max\left\{\left|\overline{P_l}\right|^2, \left|\underline{P_l}\right|^2\right\} + max\left\{\left|\overline{Q_l}\right|^2, \left|\underline{Q_l}\right|^2\right\} \leq \bar{f}_l v_{up(l)} \quad (3.e)$$

$$v_l \geq v_l^{min}, \ \hat{v}_l \leq v_l^{max} \quad (3.f)$$

$$\bar{f}_l \leq l_l^{max} \quad (3.g)$$

The MAR-OPF equations include (1) and (3) together as the set of constraints. The unknown variables are $(s^g, S, v, f, \hat{S}, \hat{v}, \underline{S}, \bar{S}, \bar{f})$. Equations (3.a) and (3.b) formulate Kirchhoff's current and voltage laws, respectively, while assuming that the power loss term ($Z_l f_l$) is zero. These yield $\hat{v}$, which is, according to Lemma I below, an upper bound for $v$. Unlike [31] and due to the presence of impedance-constant terms in (3.a), $\hat{P}$ and $\hat{Q}$ are not lower bounds for $P$ and $Q$, respectively. For deriving the lower bounds on these variables, (3.c) is introduced by modifying (3.b), in such a way that the power of constant-impedance loads is calculated using $v$ instead of $\hat{v}$. For extracting upper bounds on $P$ and $Q$, (3.d) is added, which is similar to the Kirchhoff's current low in (1.a), while instead of $f$, its upper bound $\bar{f}$ is used. The auxiliary variable $\bar{f}$ is defined in (3.e) similar to (2). However, to account for the reverse current flow in the grid, the maximum of the absolute values of upper- and lower bounds on power flow variables are used in (3.e) to define $\bar{f}$. As it can be seen, all the upper and lower bounds are independent from $f$, which is a basic requirement in modifying the OPF formulation. According to the defined upper bounds on $v$ and $f$, the voltage and current limits in MAR-OPF will be accounted for using the conservative constraints given in (3.f) and (3.g), respectively.

### D. Feasibility Analysis

In this section, we demonstrate the feasibility of the solution of the optimization problem with constraints (1) and (3). First, we define the matrices $\mathbf{M}_1$, $\mathbf{M}_2$, and $\mathbf{D}$ that are used hereafter,

$$\begin{cases} \mathbf{M}_1 = 2diag(r)\mathbf{H}diag(g) \\ \mathbf{M}_2 = 2diag(x)\mathbf{H}diag(b) \\ \mathbf{C} = [\mathbf{I} - \mathbf{G}^T + \mathbf{M}_1 + \mathbf{M}_2]^{-1} \end{cases} \quad (4)$$

$$\mathbf{D} = 2\mathbf{C}diag(r)(\mathbf{H} - \mathbf{I})diag(r) + \\ 2\mathbf{C}diag(x)(\mathbf{H} - \mathbf{I})diag(x) + \mathbf{C}diag(|z|^2) \quad (5)$$

$$\begin{cases} \pi_l = \max\left(\dfrac{P^{max}, |\mathbf{H}(p - p_{max}^g)|}{v_l^{min}}\right) \\ \varrho_l = \max\left(\dfrac{Q^{max}, |\mathbf{H}(q - q_{max}^g)|}{v_l^{min}}\right) \\ \vartheta_l = \pi_l^2 + \varrho_l^2 \end{cases} \quad (6)$$

$$\mathbf{E} = 2diag(\pi)\mathbf{H}diag(r) + 2diag(\varrho)\mathbf{H}diag(x) + diag(\vartheta)\mathbf{D} \quad (6)$$

where, the vectors $r, x, z, g, b$ and the matrices $\mathbf{G}$ and $\mathbf{H}$ are defined based on notation of section II-A.

Lemma I:

1) If the following conditions are fulfilled:

$$\|\mathbf{H}^T(-\mathbf{M}_1 - \mathbf{M}_2)\| < 1 \quad (8.a)$$

$$\mathbf{D} \geq 0 \quad (88.b)$$

then $v \leq \hat{v}$, $\underline{P} \leq P$, and $\underline{Q} \leq Q$.

2) If $(s^g, S, v, f)$ satisfies (1.c), then $f \leq \bar{f}$, $P \leq \bar{P}$, $Q \leq \bar{Q}$.

3) If $(s^g, S, v, f)$ satisfies (1.c) and $(s^g, S', v', f', \underline{S}', \bar{S}', \bar{f}')$ satisfies (1.a), (1.b), (2) and (3) with $0 \leq v' \leq v$, then there exists $(\bar{f}, \bar{S}, \underline{S})$ such that $\bar{f} \leq \bar{f}'$, $\bar{P} \leq \bar{P}'$, $\bar{Q} \leq \bar{Q}'$ $\underline{P} \leq \underline{P}'$, $\underline{Q} \leq \underline{Q}'$ and $(s^g, S, v, f, \bar{f}, \bar{S}, \underline{S})$ satisfies (3.c)-(3.e).

### Lemma II:

The feasible solution space of the optimization problem subject to (1) and (3) is a subset of the feasible solution space of the original OPF.

The proof of Lemma I is provided in Appendix II. Lemma II can be proved using the items of Lemma I. According to Lemma II, the additional constraints of MAR-OPF reduce the feasible set. Therefore, the obtained solution is not necessarily globally optimal. However, the missing portion of the feasible region covers the operation zone close to the upper voltage and current limits, where it is not desired to operate the system. In section IV, an iterative algorithm is developed to improve the quality of the obtained solution, when it resides in this region.

### Theorem I:

For every feasible solution $(s, S, f, v, \hat{v}, \hat{S}, \underline{S}, \bar{S}, \bar{f})$ of MAR-OPF there exists a feasible solution $(s, S^*, f^*, v^*)$ of the original OPF with the same power injection vectors $(s^g)$ such that:

Conditions (8.a) and (8.b) hold

$$\|\mathbf{E}\| < 1 \quad (8.c)$$

$$\exists \eta < 0.5, \quad \mathbf{D}\mathbf{E} \leq \eta\mathbf{D} \quad (8.d)$$

As noted in Appendix II, condition (8.a) guarantees that matrix $\mathbf{I} - \mathbf{G}^T + \mathbf{M}_1 + \mathbf{M}_2$ is invertible. Condition (8.b) implies that matrix $\mathbf{D}$ has nonnegative entries, which is required in the proof of Lemma I. Conditions (8.c) is necessary to ensure that there exist a power flow solution for the solution of the optimization problem. Condition (8.d) implies that the voltage magnitude of all buses increase when the current magnitude decreases at one or multiple lines of the network [31].

Theorem I states that the vector of power injections $(s^g)$ that satisfies the MAR-OPF constraints, creates voltage and current magnitudes that are within their feasible regions. Actually, if $(s^g, S, f, v, \hat{v}, \hat{S}, \underline{S}, \bar{S}, \bar{f})$ is a feasible solution of MAR-OPF, then $S, f$, and $v$ are in general not the real values of the line power, line current and nodal voltage variables (since (1.c) is replaced with (2)). The real values, $(S^*, f^*, v^*)$, are obtained by solving the set of non-linear equation (1), given $s^g$ as a known parameter. For this aim, the algorithm shown in Figure. 3 is used, where $n$ and $\epsilon$ denote the iteration number and the stopping threshold, respectively. Compared to the iterative scheme proposed in [31], the terms $Hdiag(\alpha^p)v^{(n)}$ and $Hdiag(\alpha^q)v^{(n)}$ are added to the formulation of $P^{(n)}$, and $Q^{(n)}$, in order to account for the constant-impedance loads. It should be noted that the algorithm of Figure. 3 is not executed as a part of the MAR-OPF problem. It is used solely to prove Theorem I. By using Lemma I, Lemma II, and the algorithm shown in Figure. 3, we can prove Theorem I using the same arguments given in [31]. We omit the detailed proof due to space limitations.

### III. Optimally reconfigurable Droop

In this section, the developed MAR-OPF formulation is applied to the special case of optimizing the droop parameters in a radial Active Distribution Network (ADN). An ADN interconnects consumers and IBDGs together and to the

---

Input: $\omega = (s^g, \alpha^p, \alpha^q, S, f, v, \hat{v}, \hat{S}, \underline{S}, \bar{S}, \bar{f})$.

$n = 1$: $f^{(0)} \leftarrow f, v^{(0)} \leftarrow v, P^{(0)} \leftarrow P, Q^{(0)} \leftarrow Q$

$n \geq 1$:

$$f_l^{(n)} \leftarrow \frac{\left(P_l^{(n-1)}\right)^2 + \left(Q_l^{(n-1)}\right)^2}{v_{up(l)}^{(n-1)}}$$

$$P^{(n)} \leftarrow Hp + Hdiag(g)v^{(n)} + Hdiag(r)f^{(n)}$$

$$Q^{(n)} \leftarrow Hq + Hdiag(b)v^{(n)} + Hdiag(x)f^{(n)}$$

$$v^{(n)} \leftarrow \hat{v} - Df^{(n)}$$

$$Stop \ if \left| f_l^{(n)} - f_l^{(n-1)} \right| \leq \epsilon$$

Figure. 3. The iterative ad-hoc algorithm for solving non-linear power flow equations, used for proving Theorem I.

transmission (or sub-transmission) grid through a substation bus. The voltage of this substation bus is determined by the transmission grid. A simple ADN is shown in Figure. 4. The network is assumed to be balanced, so a single-phase representation is used for the three-phase electrical network. In Appendix I, it is explained how to extend the proposed formulation to the unbalanced networks. In the simple test network shown in Figure. 4, the resource at node 2 is non-dispatchable, i.e. it does not have any control over its power injections. The resources at nodes 3 and 4 are dispatchable, that is, they control their active and reactive powers based on local voltage measurements according to the linear Q-V and P-V droop characteristics, shown in Figure. 5.

For a given IBDG at node $l$, the set of parameters of the P-V droop includes droop slope ($\alpha_l^p$) and voltage reference ($V_l^{0p}$). Those of the Q-V droop includes droop slope ($\alpha_l^q$), voltage reference ($V_l^{0q}$), and reactive power reference ($Q_l^0$). $P_{l,t}^{ava}$ is the available active power of the primary resource of the IBDG at node $l$ and time $t$. Since the droop slope has a significant impact on the generation costs and dynamic behavior of the system, it is typically set based on the dynamic analysis and economic power dispatch, which are beyond the scope of this paper [8], [32]. Therefore, we assume $\alpha_l^p$ and $\alpha_l^q$ are given and we obtain the optimal values of the remaining droop parameters to guarantee optimal operation under load and generation fluctuations.

### A. Modelling of the droop control

The aim of this section is to express the power generation of each IBDG at node $l$ and time $t$, $s_{l,t}^g = p_{l,t}^g + jq_{l,t}^g$, as a linear function of the square of the measured local voltage ($v_{l,t}$) and droop parameters.

The Taylor series expansion of function $v_{l,t} = f\left(\sqrt{v_{l,t}}\right) = \sqrt{v_{l,t}}^2$ around a given point $\sqrt{v_0}$ is:

$$v_{l,t} = v_0 + (\sqrt{v_{l,t}} - \sqrt{v_0})2\sqrt{v_0} + \left(\sqrt{v_{l,t}} - \sqrt{v_0}\right)^2 \quad (7)$$

where, $\sqrt{v_{l,t}}$ and $v_{l,t}$ represent the voltage magnitude and the squared voltage magnitude at node $l$ and time $t$. In this paper,

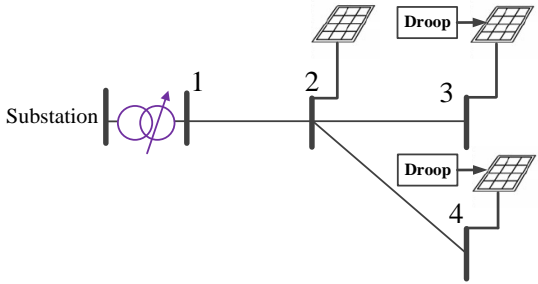

Figure. 4. A simple ADN with a Q-V and P-V droop controlled IBDG.

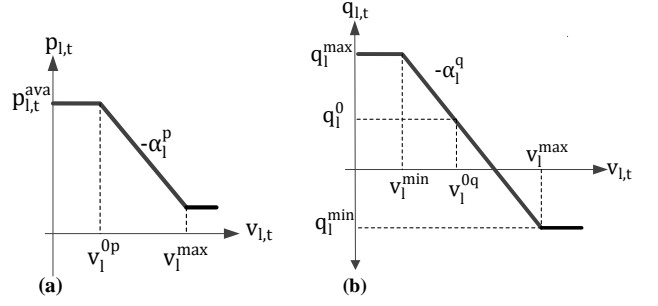

Figure. 5. The approximated P-V (a) and Q-V (b) droop characteristics.

the standard Q-V droop curve is approximated as follows:

$$q_{l,t}^g = q_{l,t}^{g0} - \alpha_l^q (v_{l,t} - v_l^{0q})$$

where the reactive power has a linear characteristic with the squared voltage magnitude (instead of the voltage magnitude). Replacing $v_{l,t}$ with its Taylor expression given in (7) results in:

$$q_{l,t}^{g0} = q_{l,t}^{g0} - \alpha_l^q \cdot 2\sqrt{v_0} \cdot \left(\sqrt{v_{l,t}} - \frac{v_l^{0q} + v_0}{2\sqrt{v_0}}\right) - \alpha_l^q \left(\sqrt{v_{l,t}} - \sqrt{v_0}\right)^2$$

If we use the first-order approximation of the Taylor series (neglecting the last expression in the right-hand side), the above formula will represent the exact Q-V droop curves, where the reactive power is proportional to the voltage magnitude. The slope of this exact droop curve is $\alpha_l^{q*} = \alpha_l^q \cdot 2\sqrt{v_0}$ and the reference voltage setpoint is $\sqrt{v_l^{0q*}} = \frac{v_l^{0q} + v_0}{2\sqrt{v_0}}$. The exact reactive power injection of the droop controller ($q_{l,t}^{g*}$) is given by:

$$q_{l,t}^{g*} = q_{l,t}^{g0} - \alpha_l^{q*} \cdot \left(\sqrt{v_{l,t}} - \sqrt{v_l^{0q*}}\right)$$

Comparing $q_{l,t}^{g*}$ with the formulation of $q_{l,t}^g$ implies that the modelled reactive power injection ($q_{l,t}^g$) is always lower than the exact one ($q_{l,t}^{g*}$) and the mismatch is:

$$q_{l,t}^g - q_{l,t}^{g*} = \alpha_l^q \left(\sqrt{v_{l,t}} - \sqrt{v_0}\right)^2$$

In other words, the modelled reactive power is a lower bound of the exact active power injection of the droop controller. The mismatch is zero if the measured voltage ($v_{l,t}$) is equal to $v_0$ (initial point of the Taylor expansion). To preserve the feasibility of the solution with respect to the maximum voltage and current limits, we choose $v_0$ equal to the maximum voltage allowed for the inverter. Therefore, we guarantee that the modelling mismatch of the power injection is zero when the voltage and/or current magnitudes are binding to their maximum limits. We apply the same analysis and solution approach also to the modeling of the P-V droop curve. As shown in [33], the approximation error is negligible in the range of active and reactive power capacities, leading to a solution with limited conservativeness. The linear approximations of the P-V and Q-V droop characteristics are shown in Figure. 5.

The P-V and Q-V droop models are given in (8), and (9), respectively. Compared to the original Q-V droop curve proposed by IEEE 1547.8 standard and depicted in Figure. 1, the considered Q-V droop curve has no dead-band. Dead-bands in Figure. 1 are usually small and we are implicitly assuming that their presence does not significantly affect the feasibility and optimality of the solution. The piecewise affine characteristic of the P-V droop is modelled using mixed-integer linear constraints (9), a procedure that is also applied in [34] for modelling hybrid systems. According to this methodology, an additional binary variable $y_{l,t}$ is assigned to each IBDG at node $l$ and time $t$. M is a sufficiently large coefficient. Constraint (9.a) states that $y_{l,t}$ equals to one when $v_{l,t}$ is larger than $v_l^{0p}$ and $y_{l,t}$ equals to zero when otherwise. If $y_{l,t} = 1$, (9.b) is imposed and (9.c) is weakened to a constraint which always holds. On the other hand, when $y_{l,t} = 0$, (9.b) is weakened and (9.c) is imposing $p_{l,t}^g$ to be equal to $p_{l,t}^{ava}$.

$$q_{l,t}^g = q_{l,t}^{g0} - \alpha_l^q (v_{l,t} - v_l^{0q}) \qquad (8)$$

$$v_l^{0p} - M(1 - y_{l,t}) \leq v_{l,t} \leq v_l^{0p} + M(y_{l,t}) \qquad (9.a)$$

$$-M(1 - y_{l,t}) \leq p_{l,t}^g - \left( p_{l,t}^{ava} - \alpha_l^p (v_{l,t} - v_l^{0p}) \right)$$
$$\leq M(1 - y_{l,t}) \qquad (9.b)$$

$$-M(y_{l,t}) \leq p_{l,t}^g - p_{l,t}^{ava} \leq M(y_{l,t}) \qquad (9.c)$$

According to the formulation above, a droop-based IBDG at bus $l$ is modelled as a constant-power generation in parallel with a constant-impedance load. The apparent power of the constant-power is equal to $(p_{l,t}^{ava} + \alpha_l^p v_l^{0p}) + j(q_{l,t}^{g0} + \alpha_l^q v_l^{0q})$, when $y_{l,t} = 1$, and equal to $p_{l,t}^{ava} + j(q_{l,t}^{g0} + \alpha_l^q v_l^{0q})$, when $y_{l,t} = 0$. The admittance of the constant-impedance load is equal to $(\alpha_l^p + j\alpha_l^q)$, when $y_{l,t} = 1$, and equal to $j\alpha_l^q$, when $y_{l,t} = 0$. Therefore, the MAR-OPF formulation developed in section II can be applied as in the following to derive the voltage and current magnitudes at each bus/line $l$ and at each time $t \in \{1,2,\ldots,T\}$.

### B. Optimal adaptive droop design

For a radial grid, we define the following optimization problem to determine the optimal droop parameters $\phi = \{V_i^{0p}, V_i^{0q}, Q_i^0 | i \in \mathcal{L}^g\}$. $\mathcal{L}^g$ denotes the set of buses hosting droop-regulated IBDGs. The set of constraints (10.e)-(10.n) are considered for each time $t$, and for each bus (or line) $l \in \mathcal{L}$ (if not mentioned otherwise).

$$\min_\phi F^{obj} = w_{pc}.F^{pc} + w_{pl}.F^{pl} + w_v.F^v \qquad (10.a)$$

$$F^{pc} = \sum_t \sum_{l \in \mathcal{L}^g} (p_{l,t}^g - p_{l,t}^{ava}) \qquad (10.b)$$

$$F^{pl} = \sum_t \sum_{l \in \mathcal{L}} r_l f_{l,t} \qquad (10.c)$$

$$F^v = \sum_{t=1}^T \sum_{l \in \mathcal{L}} |v_{l,t} - v_l^*| : \quad (|v_{l,t} - v_l^*| \geq \Delta v_l^{thr}) \qquad (10.d)$$

Subject to:

$$\text{Set of constraints (8), and (9)} \qquad \forall l \in \mathcal{L}^g \quad (10.e)$$

$$\begin{cases} v_{l,t}^{dev} \geq 0 \\ v_{l,t}^{dev} \geq v_{l,t} - (v_l^* + \Delta v_l^{thr}) \\ v_{l,t}^{dev} \geq (v_l^* - \Delta v_l^{thr}) - v_{l,t} \end{cases} \qquad (10.f)$$

$$\begin{cases} \hat{q}_{l,t}^g = q_{l,t}^{g0} - \alpha_l^q (\hat{v}_{l,t} - v_l^{0q}) \\ -M(1 - y_{l,t}) \leq \hat{p}_{l,t}^g - \left( p_{l,t}^{ava} - \alpha_l^p (\hat{v}_{l,t} - v_l^{0p}) \right) \\ \qquad\qquad \leq M(1 - y_{l,t}) \\ -M(y_{l,t}) \leq \hat{p}_{l,t}^g - p_{l,t}^{ava} \leq M(y_{l,t}) \end{cases} \qquad (10.g)$$

$$\begin{cases} s_{l,t}^g = p_{l,t}^g + jq_{l,t}^g \\ S_{l,t} = s_{l,t} - s_{l,t}^g + \sum_{m \in \mathcal{L}} G_{l,m} S_{m,t} + Z_l f_{l,t} \\ v_{l,t} = v_{up(l),t} - 2\Re(z_l^* S_{l,t}) + |z_l|^2 f_{l,t} \\ f_{l,t} \geq \dfrac{|S_{l,t}|^2}{v_{l,t}} \end{cases} \qquad (10.h)$$

$$\begin{cases} \hat{s}_{l,t}^g = \hat{p}_{l,t}^g + j\hat{q}_{l,t}^g \\ \hat{S}_{l,t} = s_{l,t} - \hat{s}_{l,t}^g + \sum_{m \in \mathcal{L}} G_{l,m} \hat{S}_{m,t} \\ \hat{v}_{l,t} = \hat{v}_{up(l),t} - 2\Re(z_l^* \hat{S}_{l,t}) \end{cases} \qquad (10.i)$$

$$\begin{cases} \overline{S_{l,t}} = s_{l,t} - s_{l,t}^g + \sum_{m \in \mathcal{L}} G_{l,m} \overline{S_{m,t}} + Z_l \overline{f_{l,t}} \\ \underline{S_{l,t}} = s_{l,t} - s_{l,t}^g + \sum_{m \in \mathcal{L}} G_{l,m} S_{m,t} \\ \max\{|\overline{P_{l,t}}|^2, |\underline{P_{l,t}}|^2\} + \max\{|\overline{Q_{l,t}}|^2, |Q_{l,t}|^2\} \leq \overline{f_{l,t}} v_{up(l),t} \end{cases} \qquad (10.j)$$

$$v_{l,t} \geq v_l^{min}, \ \hat{v}_{l,t} \leq v_l^{max} \qquad (10.k)$$

$$\overline{f_{l,t}} \leq I_l^{max} \qquad (10.l)$$

$$\begin{cases} P_{l,t} \leq \overline{P_{l,t}} \leq P_l^{max} \\ Q_{l,t} \leq \overline{Q_{l,t}} \leq Q_l^{max} \end{cases} \qquad (10.m)$$

$$\begin{cases} p_{l,t}^g \leq p_{l,max}^g \\ q_{l,min}^g \leq q_{l,t}^g \leq q_{l,max}^g \\ \dfrac{-1}{\mu_{min}}.p_{l,t}^g \leq q_{l,t}^g \leq \dfrac{1}{\mu_{min}}.p_{l,t}^g \\ (p_{l,t}^g)^2 + (q_{l,t}^g)^2 \leq (S_{l,max}^g)^2 \end{cases} \qquad (10.n)$$

The objective function in (10.a) comprises power curtailment ($F^{pc}$), power loss ($F^{pl}$), and voltage deviation ($F^v$) terms, prioritized according to the weights, $w_{pc}$, $w_{pl}$, and $w_v$, respectively. Equation (10.b) describes the total active power curtailment of each dispatchable IBDG at node $l$ and at time $t$ during overvoltage conditions. The total active power loss is described in (10.c). Equation (10.d) formulates the total deviation of the squared voltage magnitude with respect to a pre-defined target value ($v_l^*$), if beyond a certain threshold ($\Delta v_l^{thr}$). In order to obtain a convex formulation, a set of auxiliary variables $v_{l,t}^{dev}$ are introduced, which are constrained in (10.f) [35].

The security constraints of the network are given in (10.k)-(10.m), whereas (10.n) represents the power capacity limits of IBDGs. The first and second expressions of (10.m) model the active and reactive power limits that are related to the hardware of IBDG. The third expression implements the minimum allowed power factor ($\mu_{min}$) for the IBDG. Finally, the fourth expression formulates the maximum apparent power capacity ($S_{l,max}^g$) of the IBDG. This is a cone constraint exactly as the last expressions of (10.h) and (10.j).In the following, we prove that the integration of the mixed-integer linear constraints (9) into the developed MAR-OP formulation results in a feasible droop parameter set.

*Theorem II:* Under the set of conditions (8), for every feasible solution $(s^g, S, f, v, \hat{s}^g, \hat{S}, \hat{v}, \underline{S}, \bar{S}, \bar{f})$ of (10) there exists a feasible solution $(s, s^{g*}, S^*, f^*, v^*)$ of the original OPF with the same power injection vector $(s^g)$ and the same droop parameters $(\alpha^p, v^{0p}, \alpha^q, v^{0q}, q^{0g})$.

**Proof:** Let $\omega = (s^g, S, f, v, \hat{s}^g, \hat{S}, \hat{v}, \underline{S}, \bar{S}, \bar{f})$ is a feasible solution of (10). For a given IBDG at bus $l$ and time $t$, if $y_{l,t} = 1$, it means that $v_{l,t} \geq v_l^{0p}$ and $p_{l,t}^g = p_{l,t}^{ava} - \alpha_l^p(v_{l,t} - v_l^{0p})$. We adopt the ad-hoc Algorithm 1, given in section II.C., to find the corresponding power-flow solution $\omega^* = (S^*, f^*, v^*)$. According to [31], $v_{l,t}^* \geq v_{l,t}$. Therefore, $v_{l,t}^* \geq v_l^{0p}$, which follows that $p_{l,t}^{g*} = p_{l,t}^{ava} - \alpha_l^p(v_{l,t}^* - v_l^{0p})$. Therefore, the amount of power injection and load admittance in the MAR-OPF formulation used in (10), and in the power flow problem are the same. Thus, the result of Theorem I can be applied, which states that $\omega^*$ meets the security constraints.

Now, if $y_{l,t} = 0$, then $v_{l,t} \leq v_l^{0p}$ and $p_{l,t}^g = \hat{p}_{l,t}^g = p_{l,t}^{ava}$. When solving the power-flow problem with $(s, p^{ava}, \alpha^p, v^{0p}, \alpha^q, v^{0q}, q^{0g})$ the resulting solution, say $\omega^* = (S^*, f^*, v^*)$ falls into one of the following conditions:

a) $v_{l,t}^* \leq v_l^{0p}$, which yields that $p_{l,t}^{g*} = p_{l,t}^{ava}$. It means that the power injections in both MAR-OPF and power flow problems are the same. Therefore, according to Theorem I, the solution $\omega^*$ is feasible.

b) $v_{l,t}^* \geq v_l^{0p}$, which yields that $p_{l,t}^{g*} = p_{l,t}^{ava} - \alpha_l^p(v_{l,t}^* - v_l^{0p})$. If active and reactive power losses are neglected, for a given pair of buses $(m, n) \in \mathcal{L}^2$, we have $\frac{\partial \hat{v}_m}{\partial p_n^g} \geq 0$ [24]. Therefore, since $p_{l,t}^{g*} \leq \hat{p}_{l,t}^g$, $v_{l,t}^* \leq v_{l,t}$. Since $v_{l,t} \leq v_l^{max}$, we have $v_{l,t}^* \leq v_l^{max}$. Moreover, since $v_{l,t} \leq v_l^{0p} \leq v_{l,t}^*$, and since $v_l^{min} \leq v_{l,t}$, we also have $v_l^{min} \leq v_{l,t}^*$. Therefore, the solution $\omega^*$ meets the security constraints and is feasible. ∎

## IV. Optimality of the solution

As stated in section II.C, the solution of MAR-OPF does not guarantee the global optimality. In this section, we aim to improve the quality of the obtained solution by modifying (3.e), so that it is less conservative. Let define sets $\mathcal{W}^p$ and $\mathcal{W}^q$ in (11).

$$\begin{cases} \mathcal{W}^p = \left\{ (l,t) | \sum_{k \in \mathcal{L}} \left( H(diag(r) - diag(\alpha^p)D) \right)_{l,k} f_k \geq 0, |\hat{P}_l| \leq \left| \underline{P_l} \right| \right\} \\ \mathcal{W}^q = \left\{ (l,t) | \sum_{k \in \mathcal{L}} \left( H(diag(x) - diag(\alpha^q)D) \right)_{l,k} f_k \geq 0, |\hat{Q}_l| \leq \left| \underline{Q_l} \right| \right\} \end{cases} \quad (11)$$

According to (20), we conclude that $\hat{P}_{l,t} \leq P_{l,t}, \forall (l,t) \in \mathcal{W}^p$, and $\hat{Q}_{l,t} \leq Q_{l,t}, \forall (l,t) \in \mathcal{W}^q$. In, (3.e), we replace $\left| \underline{P_l} \right|$ with $|\hat{P}_l|$ for $\forall (l,t) \in \mathcal{W}^p$ and we replace $\left| \underline{Q_l} \right|$ with $|\hat{Q}_l|$ for $\forall (l,t) \in \mathcal{W}^q$. Therefore, the value of the auxiliary variable $\overline{f_{l,t}}$, that is an upper bound for the variable $f_{l,t}$, is reduced. It follows that (3.e) is less conservative and the quality of the solution, while preserving its feasibility, will be improved.

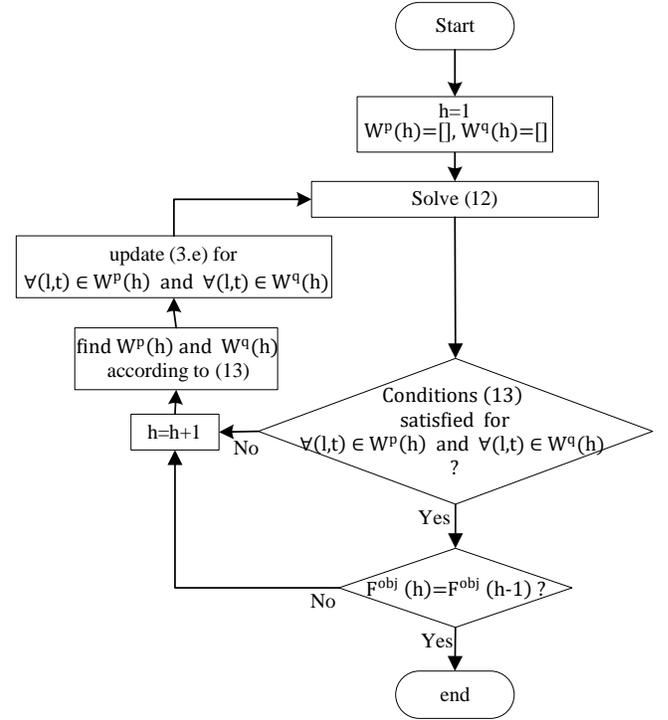

Figure. 6. Flowchart of the proposed solution algorithm for improving the solution quality (h denotes the iteration number).

To implement the proposed idea and improve the quality of the obtained solution, we propose the iterative algorithm, shown in Figure. 6. In the first iteration, the optimization problem (10) is solved. According to Theorem II, the obtained solution is feasible. Using this solution, we find $\mathcal{W}^p(1)$ and $\mathcal{W}^q(1)$ as the set of times $t$, and lines $l$, for which conditions expressed in (11) are satisfied. We update then constraint (3.e) for these pairs $(t, l)$. In the second iteration, we solve (10) with the updated constraints. Then, we check if conditions (11) are still satisfied with the obtained solution at the second iteration for $(t, l) \in \mathcal{W}^p(1)$ and $(t, l) \in \mathcal{W}^q(1)$. If all the conditions are met, then the obtained solution is feasible. Otherwise, the solution of the second iteration is not feasible and can be discarded. In both cases, we iterate the procedure by constructing $\mathcal{W}^p(2)$ and $\mathcal{W}^p(2)$ and moving to the next iteration. The iterative procedure is stopped as soon as the objective values of two consecutive iterations are the same.

## V. Simulation Results

The proposed optimization-based approach for the reconfigurable droop voltage control is tested on the 34-bus distribution network, shown in Figure. 7. This 12.66kV distribution network is a standard IEEE test benchmark, which is connected to the sub-transmission grid through a substation. The voltage of this substation is assumed fixed to 1 p.u.a., which is realized in practice using a tap changing transformer, as shown in Figure. 7. The detailed nodal and branch data is given in [36].

The original version of the IEEE 34 standard network is unbalanced. To validate the performance of the proposed formulation, we have modified the network to be balanced. First, the overall effect of mutual impedances on the three-phase system is neglected and the line impedances of all phases are assumed the same and equal to the impedance of phase A. Second, the loads at different phases of each node are assumed all equal to the load at phase A. The bus injection data is also modified by placing 11 photovoltaic (PV) resources. Eight of these PV resources are dispatchable and controlled using local P-V and Q-V droop regulators.

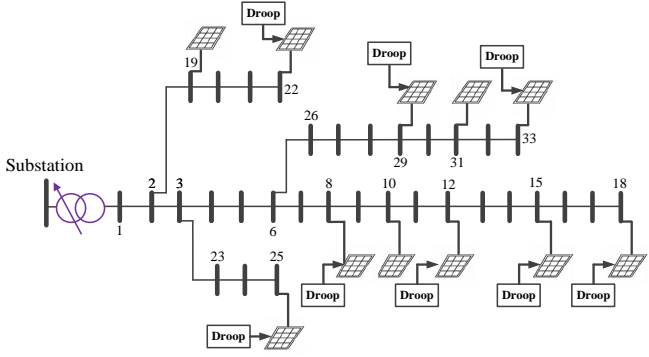

Figure. 7. One line diagram of the modified 34-bus distribution grid [36].

The time step is assumed 15 minutes. The 15-min forecast of the available active power for each PV resource on a sunny day is obtained according to the data reported in [37]. The load profile at different nodes of the network are according to the industrial, commercial, residential, and rural load patterns reported in [38]. According to ANSI C84.1 standard, the under- and overvoltage limits are set, respectively, to 0.90 and 1.05 p.u. [24]. The weighting factors of different objective terms $w_{pc}$, $w_{pl}$, and $w_v$ are set to 0.6, 0.3, and 0.1, respectively. The target voltage value ($V_t^*$) and the voltage threshold ($\Delta V^{thr}$) of all the nodes are set to 1 p.u. and 0.05 p.u., respectively.

Figure. 8 shows the nodal voltages during the whole day, when there is no control on the active and reactive powers of PV resources. It is assumed that the droop parameters can be adjusted three times a day at 7:00, 12:00, and 21:00. The optimization problem is implemented on a PC with an Intel(R) Xeon(R) CPU and 6 GB RAM; and solved in Matlab/Yalmip environment, using Gurobi solver. For handling the integrality constraints in the developed MISOCP formulation, Branch-and-Bound method is used, with the optimality gap set to 1e-10 [40].

### A. Validity of the assumptions

Note that all the conditions (8) defined in section II.C can be verified a priori for a given network. Conditions (8.a)- (8.b) are functions of the electrical parameters of the network topology, and conditions (8.c)- (8.d) are functions of the network topology

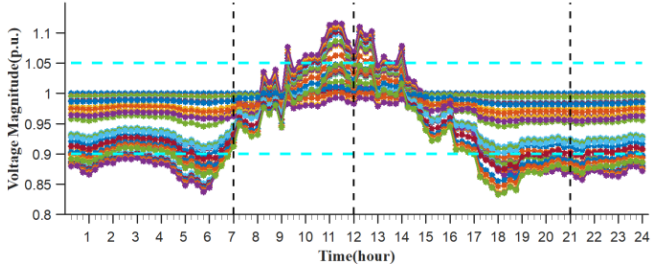

Figure. 8. 15-minute voltage profile during the day without any control on the PV resources in the 34-bus network. Dashed lines represent voltage limits. The voltages at different buses are shown in different colours.

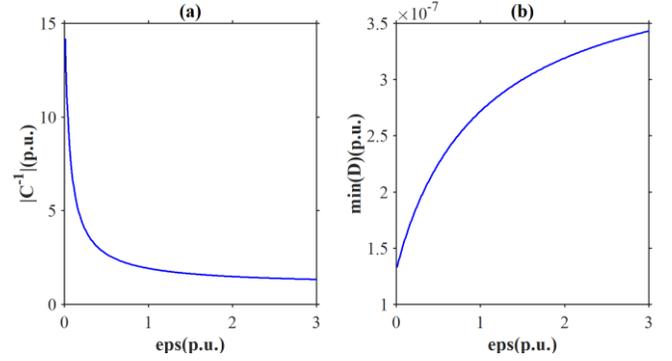

Figure. 9. Validity of conditions (8.a) (figure a) and (8.b) (figure b) for the IEEE 34-bus test network

and load/generation powers.

First, we evaluate conditions (8.a)- (8.b) for the case of the IEEE 34-bus standard test network. As proved in [41], the existence of an equilibrium operation point in the network with multiple droop-based IBDGs necessitates that $\alpha_l^p \leq \frac{1}{\sum_{k \in \mathcal{L}} \mathbf{R}_{k,l}}$ and $\alpha_l^q \leq \frac{1}{\sum_{k \in \mathcal{L}} \mathbf{X}_{k,l}}$ for all $l \in \mathcal{L}^g$, where, $\mathbf{R}_{k,l} = \sum_{h \in \mathcal{L}} r_h \mathbf{H}_{h,k} \mathbf{H}_{h,l}$ and $\mathbf{X}_{k,l} = \sum_{h \in \mathcal{L}} x_h \mathbf{H}_{h,k} \mathbf{H}_{h,l}$. Therefore, we tune the droop coefficient as follows:

$$\begin{cases} \alpha_l^p = \dfrac{1}{\sum_{k \in \mathcal{L}} \mathbf{R}_{k,l} + \epsilon} \\ \alpha_l^q = \dfrac{1}{\sum_{k \in \mathcal{L}} \mathbf{X}_{k,l} + \epsilon} \end{cases} \qquad (12)$$

where, $\epsilon \geq 0$ is tuned to make a tradeoff between convergence speed and the voltage deviation. A larger $\epsilon$ implies a faster convergence and a smaller $\epsilon$ leads to a larger voltage deviation may occur at the equilibrium. Figure. 9.a and Figure. 9.b depict, respectively, the determinant of matrix $\mathbf{C^{-1}}$ and the minimum element of matrix $\mathbf{D}$ as functions of $\epsilon$. Being nonzero and nonnegative, Figure. 9 ensures that conditions (8.a) and (8.b) hold for all possible values of $\alpha_l^p$ and $\alpha_l^q$.

To assess the validity of conditions (8.c)- (8.d), we increase the active power injections at each bus proportionally to their load share. According to this increase, the values of $p_{l,max}^g$, $q_{l,max}^g$, $P_l^{max}$, and $P_l^{max}$ are also increased. The total PV capacity is assumed 87.86 % of the total load powers. The maximum limits on the active ($P_l^{max}$) and reactive ($Q_l^{max}$) power flow of each line $l$ are assumed to be 110% of the total active and reactive loads at the downstream nodes of line $l$, respectively. The first condition, that is violated, is (8.d). However, it occurs for a total net injection equal to 5.76 MW. For this operating point, the nodal voltage-magnitudes reach a maximum value of 1.073 p.u., which is far above the maximum voltage limit.

TABLE I. Numerical results of scenario I and II using MAR-OPF and R-OPF in the 34-bus network.

| | Method | $F^{obj}$ (p.u.) | $F^{pc}$ (p.u.) | $\max_{i,t} V_{i,t}$ (p.u.) | $\min_{i,t} V_{i,t}$ (p.u.) | Computation Time (sec) |
|---|---|---|---|---|---|---|
| Scenario 1 | MAR-OPF | 0.066 | 0.016 | 1.045 | 0.9705 | 114.8 |
| | R-OPF | 0.026 | 0 | 1.079 | 0.9767 | 1.810 |
| Scenario 2 | MAR-OPF | 0.144 | 0.011 | 1.049 | 0.984 | 329.9 |
| | R-OPF | 0.245 | 0.011 | 1.108 | 0.940 | 2.594 |

### B. Comparison of MAR-OPF with R-OPF

To evaluate its performance under different network conditions, the proposed method is applied and tested for two scenarios. In scenario 1, the planning time horizon is from 07:00 to 12:00, when the PV panels are exposed to full solar irradiation and the network operates in overvoltage conditions (see Figure. 8). Scenario 2 plans for the period from 12:00 to 21:00, when the available generation of PVs vary from full capacity to zero. The aim in this scenario is to obtain a unique set of droop parameters to support both overvoltage and undervoltage conditions that are shown in Figure. 8.

The numerical results of scenarios 1 and 2 are reported in TABLE I. These results include the obtained values of the total and of the first objective terms, the computation time, and the minimum and maximum magnitudes of the voltage profile over all the nodes and over all the time steps during the planning horizon. This voltage profile is derived using power flow simulations in Matlab/MATPOWER toolbox. The simulations are run for all the time steps on the model of the network with the droop parameters that are obtained form the solution of the optimization problem. As noted in section III.A, we approximate the P-V and Q-V droop curves such that the active and reactive power have linear characteristics with the squared voltage magnitude. In order to validate among others, the accuracy of this approximation, the exact characteristic of droop curves are used in a posteriori power flow simulations. The slope and the reference setpoints of the exact droop curves are extracted from the solution of the optimization problem as explained in section III.A. The results show that the minimum (0.9 p.u.) and the maximum (1.05 p.u.) voltage limits are respected in the solution obtained using MAR-OPF approach. The short computation time for the simulation in both scenarios complies with the requirements of the operational planning in real distribution networks.

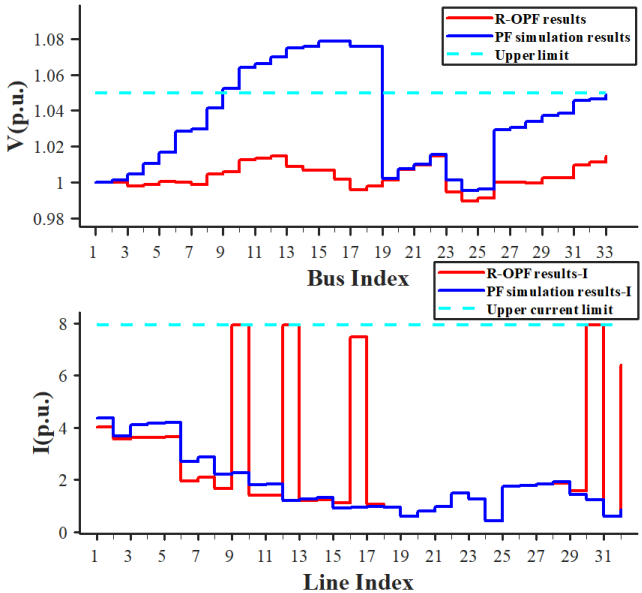

Figure. 10. The maximum voltage and current profiles in the 34-bus network under scenario 1 obtained using R-OPF. The red and blue curves represent, respectively, the values obtained from the optimization problem and the exact values obtained using a posteriori power flow simulation.

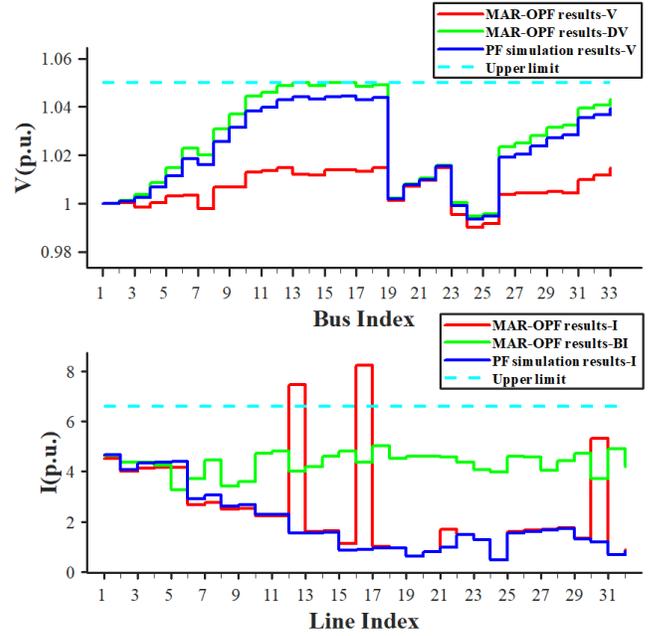

Figure. 11. The maximum voltage and current profiles in the 34-bus network under scenario 1 obtained using MAR-OPF. The red and green curves represent, respectively, the obtained values from the optimization problem for $v_{l,t}$ and $\bar{v}_{l,t}$ in the top figure, and $f_{l,t}$ and $\bar{f}_{l,t}$ in the below figure. The blue curves in the top and below figures represent, respectively, the exact values of squared voltage and current magnitudes obtained using a posteriori power flow simulation.

The performance of the MAR-OPF and R-OPF approaches are compared for the two scenarios and the results are reported in TABLE I. Although the objective value obtained with R-OPF is lower than the one of MAR-OPF, the resulting voltage profile of R-OPF approach violates the maximum voltage limit. To better understand the reason for this infeasibility, the voltage and current profiles are shown in Figure. 10. It shows the maximum values of the voltage and current variables over the planning window in scenario 1 obtained using post power flow simulations. These profiles are compared with the corresponding values obtained with the R-OPF formulation in scenario 1. As shown, there is a mismatch between these profiles, especially at the leaf nodes and lines. This is due to the reverse power flow injected by IBDGs at the leaf nodes that causes the upper voltage limit to bind and, consequently, leading to an inexact solution.

### C. Deviation from the global optimum.

In this section, we show that the compression of the feasible solution in MAR-OPF with respect to the original OPF is small. Figure. 11 shows the maximum voltage and maximum current profiles over the planning window obtained using the MAR-OPF formulation in scenario 1. As it can be seen, only the upper voltage limit is binding. In this scenario, variable $Dv$ is binding to the upper voltage limit (1.05 p.u.) at nodes 18 and 33. At these nodes, the corresponding voltage values obtained from post power flow simulation are 1.0482 and 1.0486 p.u., showing 0.17% and 0.13% of mismatch, respectively. Therefore, the solution space regarding the upper voltage limit is shrinking in the MAR-OPF just very little.

In order to assess the global optimality of the solution when line ampacity limits are binding, we run the optimization problem for t=11:00 A.M, but we relax the upper voltage limits and we increase the nodal power injections until that the first line current flow is binding to its ampacity limit. This occurs for line 1, with ampacity limit equal to 6.6 p.u.. The value of the auxiliary variable for the current flow ($\bar{f}_l$) in line 1 is equal to 6.5695 p.u., whereas the actual current value in this line is equal to 6.3212 p.u., which is equivalent to 3.7784% mismatch. It shows a very little shrinking of the feasible solution space with regarding to the ampacity limit.

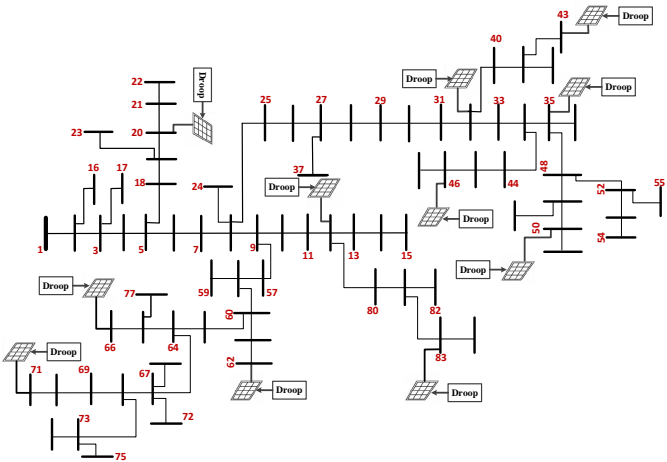

Figure. 12. One line diagram of the modified 85-bus distribution grid [42]

### D. IEEE 85-bus network

To evaluate the scalability of the proposed optimization formulation, we apply it to the MATPOWER 85-bus standard test network. The original data is reported in [42]. This network is modified in this paper to make it balanced. Moreover, controllable and non-controllable IBDGs are placed at the buses shown in Figure. 12. The characteristics of these IBDGs are the same as the ones used for the IEEE 34-bus network. We run the optimization problem for this network from 07:00 to 12:00 as in scenario 1. The simulation conditions, solar radiation data, and setting of the optimization parameters are the same as those for the 34-bus network.

The numerical results and the computation times are reported in TABLE II. As it can be seen, the proposed optimization algorithm finds the solution in a short time also in this larger network. This proves that our method can be applied to practical distribution networks with realistic sizes and complexities. The maximum voltage and currents obtained from the optimization problem over the planning horizon are depicted in Figure. 13. This figure also shows the exact values of squared voltage and current magnitudes, obtained using *a posteriori* power flow simulation. As it can be seen, both voltage and current profiles are below the maximum limits but very close to the maximum limits. The results of post power flow simulation show 1.17% and 11.12% margins with respect to the maximum allowed squared voltage and current magnitudes, respectively. These maximum values occur both at time 11:15 a.m. The maximum voltage is at node 50 and the maximum current magnitude is at line 1. These small margins indicate that the feasible solution shrinks very little with respect to the original OPF problem.

TABLE II. Numerical results of scenario 1 in the 85-bus network

| Scenario | Method | $F^{obj}$ (p.u.) | $F^{pc}$ (p.u.) | $\max\limits_{i,t} V_{i,t}$ (p.u.) | $\min\limits_{i,t} V_{i,t}$ (p.u.) | Computation Time (sec) |
|---|---|---|---|---|---|---|
| | MAR-OPF | 0.3440 | 0.0438 | 1.045 | 0.9588 | 366.21 |

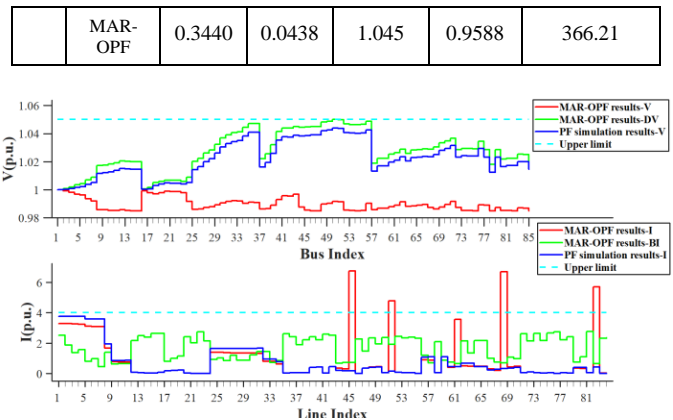

Figure. 13. The maximum voltage and current profiles in the 85-bus network under scenario 1 obtained using MAR-OPF. The red and green curves represent, respectively, the obtained values from the optimization problem for $v_{i,t}$ and $\bar{v}_{i,t}$ in the top figure, and $f_{l,t}$ and $\bar{f}_{l,t}$ in the below figure. The blue curves in the top and below figures depict, respectively, the exact values of squared voltage and current magnitudes obtained using a posteriori power flow simulation.

## VI. Conclusion

In this paper, a mixed-integer convex optimization approach is proposed for optimally designing the parameters of reconfigurable voltage droop regulators in an ADN under changing operating conditions. In this regard, first, a convex relaxation method, called MAR-OPF, is developed for the OPF problem. The feasibility of the solution under verifiable conditions is proven. In this formulation, the constant-impedance loads, upper voltage limits, and upper current limits are considered. The proposed optimization approach does not guarantee to obtain the global optimal solution, especially when it resides close to the upper current limit. For such cases, an algorithm is developed to improve the quality of the obtained solution through iterations. As a special application to the proposed relaxation method, the problem of finding the optimal parameters of P-V and Q-V droop parameters is studied. These characteristics are modelled linearly as constant-power injections in parallel with constant-impedance loads.

The contributions of the paper are validated through several simulations on a standard IEEE 34-bus test grid. First, it is numerically shown that the sufficient conditions developed for the feasibility of the MAR-OPF solution hold, with large margins, for the test network. Then, the developed optimization formulation is applied on the test network to set the parameters of the P-V and Q-V droop regulators in the network. For this aim, two scenarios are studied. In the first scenario, only the upper voltage, while in the second scenario, both upper and lower voltage limits are binding. The feasibility of the obtained solution in these scenarios is validated using power flow simulations on the test network with the obtained droop parameters from the optimization problem. Afterwards, it is illustrated that the general R-OPF formulation results in an infeasible and inexact solution when applied to the considered simulation scenarios. Finally, the distance of the obtained solutions from the global optimal solution is quantified with regard to both upper voltage limit and upper current limit. It is shown that the obtained solution is very close to global optimal solution. Therefore, the developed formulation can be used even in applications, where the global optimality of the solution is concerned.

## VII. Appendix I: Extension of the optimization model to unbalanced grids

As noted in Section III, the optimization problem (10) is formulated for balanced networks. However, distribution networks are often unbalanced. In this appendix, we explain how our method can be extended to this case. We still assume that i) all the impedances of a given line at different phases are equal, and ii) the mutual impedances between different phases of each line are identical. These assumptions generally hold, especially in medium-voltage distribution grids where the lines/cables of different phases are the same and placed in a symmetrical configuration. In this regard, all the voltage, current, and power flow variables are decomposed by using the well-known sequence transformation method. As a result, the unbalanced grid is decomposed into three symmetrical and balanced three-phase circuits. For each of these circuits, we solve problem (10).

More in details, we transform the voltage/current limits from the phase domain to the sequence domain according to the methodology given in [31]. We assume that the negative and zero sequence terms of voltage and current magnitudes are binding to their standard/normal limits. Therefore, the voltage and current limits associated with the positive sequence terms are derived *a priori*. The network security constraints ((10.l) and (10.m)) related to three sequences are integrated into the optimization problem. Once the optimization problem is solved, we transform the obtained set of electrical state variables from the sequence domain back into the phase domain.

## VIII. Appendix II: Proof of Lemma I

Equations (13) and (14) are the matrix form of the line power balance (1.a) and nodal voltage (1.b) equations, respectively.

$$\begin{cases} P = \mathbf{H}p + \mathbf{H}diag(r)f + \mathbf{H}diag(g)v \\ Q = \mathbf{H}q + \mathbf{H}diag(x)f + \mathbf{H}diag(b)v \end{cases} \quad (13)$$

$$v = \mathbf{G}^T v + v_0 e - 2diag(r)P - \\ 2diag(x)Q + diag(|z|^2) \quad (14)$$

Eliminating $P$ and $Q$ from (14) using (13) yields the following:

$$(\mathbf{I} - \mathbf{G}^T + \mathbf{M}_1 + \mathbf{M}_2)v \\ = v_0 e - 2diag(r)\mathbf{H}p - 2diag(x)\mathbf{H}q \\ + (-2diag(r)\mathbf{H}diag(r) - 2diag(x)\mathbf{H}diag(x \\ + diag(|z|^2)) \quad (15)$$

where, $e = (1, 0, ..., 0)^T$. Matrices $\mathbf{M}_1$ and $\mathbf{M}_2$ are defined in (4). According to the identity theorem and under condition (8.a) it follows [42]:

$$[\mathbf{I} - \mathbf{H}^T(-\mathbf{M}_1 - \mathbf{M}_2)].[\mathbf{I} + \mathbf{H}^T(-\mathbf{M}_1 - \mathbf{M}_2) \\ + \cdots + \left(\mathbf{H}^T(-\mathbf{M}_1 - \mathbf{M}_2)\right)^n] = \mathbf{I}$$

Therefore, $[\mathbf{I} - \mathbf{H}^T(-\mathbf{M}_1 - \mathbf{M}_2)]$ is invertible. Since $\mathbf{H} = (\mathbf{I} - \mathbf{G})^{-1}$, we conclude that $\mathbf{C} = [\mathbf{I} - \mathbf{G}^T + \mathbf{M}_1 + \mathbf{M}_2]^{-1}$ exist

when (8.a) holds. Therefore, (15) can be solved in terms of $v$ as follows:

$$v = v_0 \mathbf{C}e - 2\mathbf{C}diag(r)\mathbf{H}p - 2\mathbf{C}diag(x)\mathbf{H}q - \mathbf{D}f \quad (16)$$

where, the matrix $\mathbf{D}$ is defined in (8).

In the same way, (3.a) and (3.b) are rewritten as follows:

$$\begin{cases} \hat{P} = \mathbf{H}p + \mathbf{H}diag(g)\hat{v} \\ \hat{Q} = \mathbf{H}q + \mathbf{H}diag(b)\hat{v} \end{cases} \quad (17)$$

$$\hat{v} = v_0 \mathbf{C}e - 2\mathbf{C}diag(r)\mathbf{H}p - 2\mathbf{C}diag(x)\mathbf{H}q \quad (18)$$

Comparing (16) and (18) yields:

$$v = \hat{v} - \mathbf{D}f \quad (19)$$

Therefore, under conditions (8.a) and (8.b), $v \le \hat{v}$. Comparing (13) and (17) yields:

$$\begin{cases} P = \hat{P} + \mathbf{H}(diag(r) - diag(g)\mathbf{D})f \\ Q = \hat{Q} + \mathbf{H}(diag(x) - diag(b)\mathbf{D})f \end{cases} \quad (20)$$

Using (3.c) and (3.d), we have:

$$\begin{cases} \underline{P} = P - \mathbf{H}diag(r)f \\ \underline{Q} = Q - \mathbf{H}diag(x)f \end{cases} \quad (21)$$

$$\begin{cases} \bar{P} = P + \mathbf{H}diag(r)(\bar{f} - f) \\ \bar{Q} = Q + \mathbf{H}diag(x)(\bar{f} - f) \end{cases} \quad (22)$$

According to (21), $\underline{P} \le P$ and $\underline{Q} \le Q$. Therefore, item 1 of Lemma I holds.

We prove items 2 and 3 by induction, starting from the leaves of the grid. Let define $height(l)$ as equal to zero if $l$ is a leaf node and otherwise equal to $1 + max(height(up(l)))$.

**Base case (height = 0):** For a leaf node $l$, we show that Lemma I holds.

Item 2: Using (3.c) and (3.d) for a leaf node $l$, we have:

$$\begin{cases} \underline{P}_l = P_l - r_l f_l = \bar{P}_l - r_l \bar{f}_l = p_l + v_l g_l \\ \underline{Q}_l = Q_l - x_l f_l = \bar{Q}_l - x_l \bar{f}_l = q_l + v_l b_l \end{cases} \quad (23)$$

According to (3.e):

$$\bar{f}_l \ge \frac{(|p_l + v_l g_l|^2 + |q_l + v_l b_l|^2)}{v_l}$$

Since $(s, S, v, f)$ satisfies (1.c):

$$f_l = \frac{(P_l - r_l f_l)^2 + (Q_l - x_l f_l)^2}{v_l} = \\ \frac{(p_l + v_l g_l)^2 + (q_l + v_l b_l)^2}{v_l}$$

Therefore, for a leaf node $l$, $\bar{f}_l \ge f_l$ and according to (22), $P_l \le \bar{P}_l$ and $Q_l \le \bar{Q}_l$.

Item 3: since $v'_l \le v_l$, according to (3.d), it follows that $\underline{P}'_l \le \underline{P}_l$ and $\underline{Q}'_l \le \underline{Q}_l$. since $(\underline{S}', \bar{S}', \bar{f}')$ satisfies (3.e), using (23), we have:

$$\overline{f'_l} \ge \frac{|p_l + v'_l g_l|^2 + |q_l + v'_l b_l|^2}{v'_l} \quad (24)$$

Let us choose $\bar{f}_l$ as follows:

$$\overline{f}_l = \frac{|p_l + v_l g_l|^2 + |q_l + v_l b_l|^2}{v_l} \quad (25)$$

According to (23), $\overline{f}_l$ satisfies (3.e). Now, we prove that $\overline{f}_l \leq \overline{f}_l'$. Since $v_l' \leq v_l$, it follows that $-p_l - v_l g_l \leq -p_l - v_l' g_l$. It is logical to assume that an IBDG at a leaf node $l$ injects a positive net power to the grid. It means that $|p_l + v_l \alpha_l^p| = -p_l - v_l g_l$. Therefore, $|p_l + v_l g_l| \leq |p_l + v_l' g_l|$. The same applies also to the reactive power. Thus, according to (24) and (25), we conclude that $\overline{f}_l \leq \overline{f}_l'$ and using (23) it follows that $\overline{P}_l \leq \overline{P}_l'$ and $\overline{Q}_l \leq \overline{Q}_l'$.

**Induction Step:** Assume that the statements of Lemma I are true for all the buses with height $\leq n$. Let $k$ be the bus with $height = n + 1$.

Item 2: For all the downstream buses $l$, we have $height(l) \leq n$. Therefore, $P_l \leq \overline{P}_l$ and $Q_l \leq \overline{Q}_l$. It follows that

$$\underline{P}_k \leq P_k - r_k f_k \leq \overline{P}_k - r_k \overline{f}_k$$

And

$$\underline{Q}_k \leq Q_k - x_k f_k \leq \overline{Q}_k - x_k \overline{f}_k.$$

Thus,

$$|P_k - r_k f_k|^2 + |Q_k - x_k f_k|^2 \leq$$
$$\max \left\{ |\overline{P}_k - r_k \overline{f}_k|^2, |\underline{P}_k|^2 \right\} + \max \left\{ |\overline{Q}_k - x_k \overline{f}_k|^2, |\underline{Q}_k|^2 \right\}$$

According to (3.e),

$$\overline{f}_k v_k \geq |P_k - r_k f_k|^2 + |Q_k - x_k f_k|^2 = f_k v_k.$$

According to (22), it yields $P_k \leq \overline{P}_k$ and $Q_k \leq \overline{Q}_k$. Therefore, item 2 of Lemma I is proved in both basis and induction steps.

Item 3: since $v_k' \leq v_k$, according to (3.d), it follows that $\underline{P}_k' \leq \underline{P}_k$ and $\underline{Q}_k' \leq \underline{Q}_k$. Using (3.e) for the line $k$, we have:

$$\overline{f}_k \geq \frac{\max \left\{ |\overline{P}_k' - r_k \overline{f}_k'|, |\underline{P}_k'| \right\}^2 + \max \left\{ |\overline{Q}_k' - x_k \overline{f}_k'|, |\underline{Q}_k'| \right\}^2}{v_k'} \quad (26)$$

Let us choose $\overline{f}_k$ as follows:

$$\overline{f}_k = \frac{\max \left\{ |\overline{P}_k - r_k \overline{f}_k|, |\underline{P}_k| \right\}^2 + \max \left\{ |\overline{Q}_k - r_k \overline{f}_k|, |\underline{Q}_k| \right\}^2}{v_k} \quad (27)$$

Therefore, $\overline{f}_k$ satisfies (3.e). Now, we prove that $\overline{f}_k \leq \overline{f}_k'$. Using (3.c) and (3.d), we rewrite the numerator of the expression at the right hand side of (27) as follows:

$$\max \left\{ \left| \sum_{l \in \mathcal{L}} G_{k,l} \overline{P}_l \right|, \left| \sum_{l \in \mathcal{L}} G_{k,l} \underline{P}_l \right| \right\}^2 +$$
$$\max \left\{ \left| \sum_{l \in \mathcal{L}} G_{k,l} \overline{Q}_l \right|, \left| \sum_{l \in \mathcal{L}} G_{k,l} \underline{Q}_l \right| \right\}^2 +$$
$$\left[ 2(p_k + g_k v_k) \sum_{l \in \mathcal{L}} G_{k,l} \overline{P}_l + 2(q_k + b_k v_k) \sum_{l \in \mathcal{L}} G_{k,l} \overline{Q}_l \right] +$$
$$\left[ (p_k + g_k v_k)^2 + (q_k + b_k v_k)^2 \right]$$

In the following, we evaluate the three terms, given above:

- According to the induction assumption for all the downstream lines $l$, we have $\overline{P}_l' \geq \overline{P}_l$, and $\underline{P}_l' \leq \underline{P}_l$. Therefore,

$$\left( \frac{1}{v_k} \right) \max \left\{ \left| \sum_{l \in \mathcal{L}} G_{k,l} \overline{P}_l \right|, \left| \sum_{l \in \mathcal{L}} G_{k,l} \underline{P}_l \right| \right\} \leq$$
$$\left( \frac{1}{v_k'} \right) \max \left\{ \left| \sum_{l \in \mathcal{L}} G_{k,l} \overline{P}_l' \right|, \left| \sum_{l \in \mathcal{L}} G_{k,l} \underline{P}_l' \right| \right\}$$

Note that $v_k' \leq v_k$. The same result applies also to the reactive power flow.

- Applying the similar argument, it follows that:

$$2 \left( \frac{p_k}{v_k} + g_k \right) \sum_{l \in \mathcal{L}} G_{k,l} \overline{P}_l \leq 2 \left( \frac{p_k'}{v_k'} + g_k \right) \sum_{l \in \mathcal{L}} G_{k,l} \overline{P}_l'$$

The same result also applies to the reactive power flow.

- As we proved in the basic step for a given bus $k$,

$$\frac{(p_k + g_k v_k)^2 + (q_k + b_k v_k)^2}{v_k} \leq$$
$$\frac{(p_k + g_k v_k')^2 + (q_k + b_k v_k')^2}{v_k'}$$

Therefore, $\overline{f}_k' \geq \overline{f}_k$. Consequently, from (3.c), we have $\overline{P}_k \leq \overline{P}_k'$ and $\overline{Q}_k \leq \overline{Q}_k'$. Both basis and induction steps are proved, which completes the proof of item 3 of Lemma I.

## IX. References

[1] H. Sekhavatmanesh und R. Cherkaoui, „Distribution Network Restoration in a Multiagent Framework Using a Convex OPF Model", *IEEE Transactions on Smart Grid*, Bd. 10, Nr. 3, S. 2618–2628, Mai 2019, doi: 10.1109/TSG.2018.2805922.

[2] E. Dall'Anese, S.V. Dhople, und G. B. Giannakis, „Optimal Dispatch of Photovoltaic Inverters in Residential Distribution Systems", *IEEE Transactions on Sustainable Energy*, Bd. 5, Nr. 2, S. 487–497, Apr. 2014, doi: 10.1109/TSTE.2013.2292828.

[3] H. J. Liu, W. Shi, und H. Zhu, „Decentralized Dynamic Optimization for Power Network Voltage Control", *IEEE Transactions on Signal and Information Processing over Networks*, Bd. 3, Nr. 3, S. 568–579, Sep. 2017, doi: 10.1109/TSIPN.2016.2631886.

[4] H. Zhu und H. J. Liu, „Fast Local Voltage Control Under Limited Reactive Power: Optimality and Stability Analysis", *IEEE Transactions on Power Systems*, Bd. 31, Nr. 5, S. 3794–3803, Sep. 2016, doi: 10.1109/TPWRS.2015.2504419.

[5] Standards Coordinating Committee 21 of Institute of Electrical and Electronics Engineers, „IEEE Draft Recommended Practice for Establishing Methods and Procedures that Provide Supplemental Support for Implementation Strategies for Expanded Use of IEEE Standard 1547", IEEE Standard, P1547.8TM/D8, Aug. 2014.

[6] S. Bolognani und S. Zampieri, „A Distributed Control Strategy for Reactive Power Compensation in Smart Microgrids", *IEEE Transactions on Automatic Control*, Bd. 58, Nr. 11, S. 2818–2833, Nov. 2013, doi: 10.1109/TAC.2013.2270317.

[7] J. Zhao und F. Dörfler, „Distributed control and optimization in DC microgrids", *Automatica*, Bd. 61, S. 18–26, Nov. 2015, doi: 10.1016/j.automatica.2015.07.015.

[8] J. W. Simpson-Porco, F. Dörfler, und F. Bullo, „Synchronization and power sharing for droop-controlled inverters in islanded microgrids", *Automatica*, Bd. 49, Nr. 9, S. 2603–2611, Sep. 2013, doi: 10.1016/j.automatica.2013.05.018.

[9] Y. Zhang, X. Meng, A. M. Shotorbani, und L. Wang, „Minimization of AC-DC Grid Transmission Loss and DC Voltage Deviation Using Adaptive Droop Control and Improved AC-DC Power Flow Algorithm", *IEEE Transactions on Power Systems*, Bd. 36, Nr. 1, S. 744–756, Jan. 2021, doi: 10.1109/TPWRS.2020.3020039.

[10] S. Karagiannopoulos, P. Aristidou, und G. Hug, „Data-Driven Local Control Design for Active Distribution Grids Using Off-Line Optimal Power Flow and Machine Learning Techniques", *IEEE Transactions on Smart Grid*, Bd. 10, Nr. 6, S. 6461–6471, Nov. 2019, doi: 10.1109/TSG.2019.2905348.

[11] K. Baker, A. Bernstein, E. Dall'Anese, und C. Zhao, „Network-Cognizant Voltage Droop Control for Distribution Grids", *IEEE Transactions on Power Systems*, Bd. 33, Nr. 2, S. 2098–2108, März 2018, doi: 10.1109/TPWRS.2017.2735379.

[12] Y. Xia, Y. Peng, P. Yang, M. Yu, und W. Wei, „Distributed Coordination Control for Multiple Bidirectional Power Converters in a Hybrid AC/DC Microgrid", *IEEE Transactions on Power Electronics*, Bd. 32, Nr. 6, S. 4949–4959, Juni 2017, doi: 10.1109/TPEL.2016.2603066.

[13] Z. Liu, M. Su, Y. Sun, X. Zhang, X. Liang, und M. Zheng, „A Comprehensive Study on Existence and Stability of Equilibria of DC Distribution Networks with Constant Power Loads", *IEEE Transactions on Automatic Control*, S. 1–1, 2021, doi: 10.1109/TAC.2021.3072084.

[14] Y. Guo, D. J. Hill, und Y. Wang, „Nonlinear decentralized control of large-scale power systems", *Automatica*, Bd. 36, Nr. 9, S. 1275–1289, Sep. 2000, doi: 10.1016/S0005-1098(00)00038-8.

[15] G. Agundis-Tinajero *u. a.*, „Extended-Optimal-Power-Flow-Based Hierarchical Control for Islanded AC Microgrids", *IEEE Transactions on Power Electronics*, Bd. 34, Nr. 1, S. 840–848, Jan. 2019, doi: 10.1109/TPEL.2018.2813980.

[16] S. H. Low, „Convex Relaxation of Optimal Power Flow—Part I: Formulations and Equivalence", *IEEE Transactions on Control of Network Systems*, Bd. 1, Nr. 1, S. 15–27, März 2014, doi: 10.1109/TCNS.2014.2309732.

[17] S. Bose, S. H. Low, T. Teeraratkul, und B. Hassibi, „Equivalent Relaxations of Optimal Power Flow", *IEEE Transactions on Automatic Control*, Bd. 60, Nr. 3, S. 729–742, März 2015, doi: 10.1109/TAC.2014.2357112.

[18] R. Madani, J. Lavaei, und R. Baldick, „Convexification of Power Flow Equations in the Presence of Noisy Measurements", *IEEE Transactions on Automatic Control*, Bd. 64, Nr. 8, S. 3101–3116, Aug. 2019, doi: 10.1109/TAC.2019.2897939.

[19] S. Magnússon, P. C. Weeraddana, und C. Fischione, „A Distributed Approach for the Optimal Power-Flow Problem Based on ADMM and Sequential Convex Approximations", *IEEE Transactions on Control of Network Systems*, Bd. 2, Nr. 3, S. 238–253, Sep. 2015, doi: 10.1109/TCNS.2015.2399192.

[20] Z. Miao, L. Fan, H. G. Aghamolki, und B. Zeng, „Least Squares Estimation Based SDP Cuts for SOCP Relaxation of AC OPF", *IEEE Transactions on Automatic Control*, Bd. 63, Nr. 1, S. 241–248, Jan. 2018, doi: 10.1109/TAC.2017.2719607.

[21] J. Lavaei, D. Tse, und B. Zhang, „Geometry of Power Flows and Optimization in Distribution Networks", *IEEE Transactions on Power Systems*, Bd. 29, Nr. 2, S. 572–583, März 2014, doi: 10.1109/TPWRS.2013.2282086.

[22] L. Gan, N. Li, U. Topcu, und S. Low, „On the exactness of convex relaxation for optimal power flow in tree networks", in *2012 IEEE 51st IEEE Conference on Decision and Control (CDC)*, Dez. 2012, S. 465–471. doi: 10.1109/CDC.2012.6426045.

[23] S. Huang, Q. Wu, J. Wang, und H. Zhao, „A Sufficient Condition on Convex Relaxation of AC Optimal Power Flow in Distribution Networks", *IEEE Transactions on Power Systems*, Bd. 32, Nr. 2, S. 1359–1368, März 2017, doi: 10.1109/TPWRS.2016.2574805.

[24] L. Gan, N. Li, U. Topcu, und S. H. Low, „Exact Convex Relaxation of Optimal Power Flow in Radial Networks", *IEEE Transactions on Automatic Control*, Bd. 60, Nr. 1, S. 72–87, Jan. 2015, doi: 10.1109/TAC.2014.2332712.

[25] S. Brahma, N. Nazir, H. Ossareh, und M. R. Almassalkhi, „Optimal and Resilient Coordination of Virtual Batteries in Distribution Feeders", *IEEE Transactions on Power Systems*, Bd. 36, Nr. 4, S. 2841–2854, Juli 2021, doi: 10.1109/TPWRS.2020.3043632.

[26] N. Nazir und M. Almassalkhi, „Voltage Positioning Using Co-Optimization of Controllable Grid Assets in Radial Networks", *IEEE Transactions on Power Systems*, Bd. 36, Nr. 4, S. 2761–2770, Juli 2021, doi: 10.1109/TPWRS.2020.3044206.

[27] N. Nazir und M. Almassalkhi, „Grid-aware aggregation and realtime disaggregation of distributed energy resources in radial networks", *arXiv:1907.06709 [math]*, Okt. 2021, Zugegriffen: 20. Januar 2022. [Online]. Verfügbar unter: http://arxiv.org/abs/1907.06709

[28] D. Lee, H. D. Nguyen, K. Dvijotham, und K. Turitsyn, „Convex Restriction of Power Flow Feasibility Sets", *IEEE Transactions on Control of Network Systems*, Bd. 6, Nr. 3, S. 1235–1245, Sep. 2019, doi: 10.1109/TCNS.2019.2930896.

[29] D. Lee, K. Turitsyn, D. K. Molzahn, und L. A. Roald, „Feasible Path Identification in Optimal Power Flow With Sequential Convex

[30] D. K. Molzahn, B. C. Lesieutre, und L. A. DeMarco, „Approximate Representation of ZIP Loads in a Semidefinite Relaxation of the OPF Problem", *IEEE Transactions on Power Systems*, Bd. 29, Nr. 4, S. 1864–1865, Juli 2014, doi: 10.1109/TPWRS.2013.2295167.

[31] M. Nick, R. Cherkaoui, J. L. Boudec, und M. Paolone, „An Exact Convex Formulation of the Optimal Power Flow in Radial Distribution Networks Including Transverse Components", *IEEE Transactions on Automatic Control*, Bd. 63, Nr. 3, S. 682–697, März 2018, doi: 10.1109/TAC.2017.2722100.

[32] J. Schiffer, R. Ortega, A. Astolfi, J. Raisch, und T. Sezi, „Conditions for stability of droop-controlled inverter-based microgrids", *Automatica*, Bd. 50, Nr. 10, S. 2457–2469, Okt. 2014, doi: 10.1016/j.automatica.2014.08.009.

[33] F. Sun, J. Ma, M. Yu, und W. Wei, „A Robust Optimal Coordinated Droop Control Method for Multiple VSCs in AC–DC Distribution Network", *IEEE Transactions on Power Systems*, Bd. 34, Nr. 6, S. 5002–5011, Nov. 2019, doi: 10.1109/TPWRS.2019.2919904.

[34] A. Bemporad und M. Morari, „Control of systems integrating logic, dynamics, and constraints", *Automatica*, Bd. 35, Nr. 3, S. 407–427, März 1999, doi: 10.1016/S0005-1098(98)00178-2.

[35] H. Sekhavatmanesh und R. Cherkaoui, „Analytical Approach for Active Distribution Network Restoration Including Optimal Voltage Regulation", *IEEE Transactions on Power Systems*, Bd. 34, Nr. 3, S. 1716–1728, Mai 2019, doi: 10.1109/TPWRS.2018.2889241.

[36] M. E. Baran und F. F. Wu, „Network reconfiguration in distribution systems for loss reduction and load balancing", *IEEE Transactions on Power Delivery*, Bd. 4, Nr. 2, S. 1401–1407, Apr. 1989, doi: 10.1109/61.25627.

[37] „Solar Irradiation Data (SODA)", 15. Februar 2021. http://soda-pro.com

[38] E. Lopez, H. Opazo, L. Garcia, und P. Bastard, „Online reconfiguration considering variability demand: applications to real networks", *IEEE Transactions on Power Systems*, Bd. 19, Nr. 1, S. 549–553, Feb. 2004, doi: 10.1109/TPWRS.2003.821447.

[39] „Std, A. N. S. I. ‚C84. 1-2011.' American National Standard for Electric Power Systems and Equipment-Voltage Ratings (60 Hertz) (2011)."

[40] H. Sekhavatmanesh und R. Cherkaoui, „A Novel Decomposition Solution Approach for the Restoration Problem in Distribution Networks", *IEEE Transactions on Power Systems*, Bd. 35, Nr. 5, S. 3810–3824, Sep. 2020, doi: 10.1109/TPWRS.2020.2982502.

[41] M. Farivar, L. Chen, und S. Low, „Equilibrium and dynamics of local voltage control in distribution systems", in *52nd IEEE Conference on Decision and Control*, Dez. 2013, S. 4329–4334. doi: 10.1109/CDC.2013.6760555.

[42] P. Shiu, „Complex variables: Introduction and applications, by Mark J. Ablowitz and Athanassios S. Fokas. Pp. 647. $19.95. 1997. ISBN 0 521 48523 1 (Cambridge University Press).", *The Mathematical Gazette*, Bd. 83, Nr. 496, S. 183–184, März 1999, doi: 10.2307/3618749.

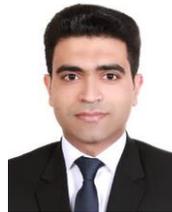

**Hossein Sekhavatmanesh** (M'14) received the PhD degree in Electrical Engineering from EPFL, Lausanne, Switzerland in 2020. He received the M.Sc. degree in power system engineering from Sharif University of Technology, Tehran, Iran in 2015 and the B.Sc. degree in electrical engineering from Amirkabir University of Technology, Tehran, Iran in 2013.

From 2013 to 2015, he was a technical expert of the R&D in MehadSanat Co. (a startup in energy), Tehran, Iran and from 2012 to 2013, he worked as a consultant engineer at TUSRC (the metro system operator), Tehran, Iran. Since October 2020, he has been a scientific collaborator with FHNW, Windisch, Switzerland.

His research interest includes optimal design and optimal operation of industrial networks, power system studies for industrial plants, and distributed control of power converters. He has been serving as a reviewer of many IEEE Transaction journals in control and power system community for more than 5 years.

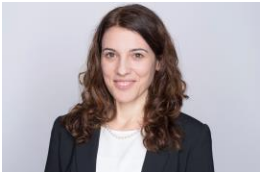

**Silvia Mastellone** is Professor for Control and Signal Processing at the University of Applied Science Northwestern Switzerland. She obtained her PhD degree in Systems and Entrepreneurial Engineering from the University of Illinois at Urbana-Champaign in 2008. From 2008 to 2016 she was Principal Scientist at ABB Corporate Research Center in Switzerland, where she led research projects in the area of advanced control for energy systems. Her research interests include nonlinear and decentralized control and estimation and networked control systems, with applications in power conversion and energy systems. She is a member of the IFAC Industry Executive Committee and a member of the advisory board for the multiutility IBB.

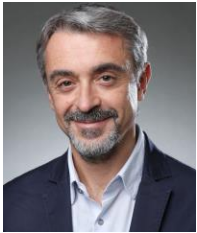

**Giancarlo Ferrari-Trecate** (SM'12) received the Ph.D. degree in Electronic and Computer Engineering from the Universita' degli Studi di Pavia in 1999. Since September 2016 he is Professor at EPFL, Lausanne, Switzerland. In spring 1998, he was a Visiting Researcher at the Neural Computing Research Group, University of Birmingham, UK. In fall 1998, he joined as a Postdoctoral Fellow the Automatic Control Laboratory, ETH, Zurich, Switzerland. He was appointed Oberassistent at ETH, in 2000. In 2002, he joined INRIA, Rocquencourt, France, as a Research Fellow. From March to October 2005, he was researcher at the Politecnico di Milano, Italy. From 2005 to August 2016, he was Associate Professor at the Dipartimento di Ingegneria Industriale e dell'Informazione of the Universita' degli Studi di Pavia.

His research interests include scalable control, microgrids, machine learning, networked control systems and hybrid systems.

Giancarlo Ferrari Trecate was the recipient of the Researcher Mobility Grant from the Italian Ministry of Education, University and Research in 2005. He is currently a member of the IFAC Technical Committees on Control Design and Optimal Control, and the Technical Committee on Systems Biology of the IEEE SMC society. He has been serving on the editorial board of Automatica for nine years and of Nonlinear Analysis: Hybrid Systems.